\documentclass[a4paper,12pt,oneside]{amsart}       
\usepackage{graphicx}

\usepackage{amssymb,amsmath,amsthm,mathrsfs}
\usepackage{svg}
\usepackage[textwidth=14cm]{geometry}

\usepackage[hidelinks]{hyperref}

\hypersetup{
    colorlinks=true,
    linkcolor=black,
    citecolor=black,
    filecolor=black,
    urlcolor=black,
}

\title[A. Aizenbud, N. Avni and R. Rubio]{Correction to: Representation  Growth\\  and Rational Singularities of the \\ Moduli Space of Local Systems}

\thanks{The third author has been funded by the Marie Sk\l odowska-Curie grant agreement No 750885 and the Spanish R\&D grant PID2019-109339GA-C32.\\ 
The third author thanks S. Carmeli for introducing the basics of \cite{aizav} to him, and the Weizmann Institute of Science for its hospitality.}

\author[A. Aizenbud]{Avraham Aizenbud}
\author[N. Avni]{Nir Avni}
\author[R. Rubio]{Roberto Rubio}

\address{A. Aizenbud \\
Weizmann Institute of Science, 76100 Rehovot, Israel }

\address{N. Avni\\ Northwestern University, Evanston, IL 60201, USA}

\address{R. Rubio \\ 
Universitat de Barcelona, 08007 Barcelona, Spain, and  \phantom{Universi-}
Universitat Aut\`onoma de Barcelona, 08193 Barcelona, Spain}

\email{roberto.rubio@ub.edu}

\begin{document}

\maketitle

\begin{abstract}
 We explain and correct a mistake in Section 2.6 and Appendix C of the first and second author's paper ``Representation  Growth  and Rational Singularities of the Moduli Space of Local Systems''\cite{aizav}.
\end{abstract}

\section*{}\vspace{-.7cm}

We use throughout the notation and conventions of \cite{aizav}. The source of the mistake is the description of the set $S_2$ in page 272. The elements
\begin{equation}\label{eq:points}
(\{ (d,d-1),(d-1,d-1)\},(d,d-1)),\; ( \{(d-1,d-1),(d-1,d)\},(d-1,d) )     
\end{equation}
are not considered. This spoils the proof of Lemma 2.40 (the blue and green subgraphs are no longer trees), which is used to prove Theorem 2.1.

\section{A straightforward correction}\label{sec:correction}

The right description of $S_2$ is
\begin{multline*}
S_2= \left\{(\{(i,j),(j,{l})\},(i,{l})) \in  I^{(2)} \times J\, |\,  \{ i, l\}=\{d-1,d\} \textrm{ and } j=  \left \lfloor \frac{i+{l}}{2} \right\rfloor,\right.\\ \left. \textrm{ or } \{ i, l\}\neq \{d-1,d\} \textrm{ and } j=  \left \lceil \frac{i+{l}}{2}\right \rceil +\delta_{i,{l}}   \right\} .
\end{multline*}
Then, $\Gamma_2$ is the polygraph\ attached to the graph $\Gamma_3=(I,E)$, with 
\begin{multline*}E= \left\{\{(i,j),(j,{l})\} \in I^{(2)} \, | \,   \{ i, l\}=\{d-1,d\} \textrm{ and } j=  \left \lfloor \frac{i+{l}}{2} \right\rfloor, \right. \\ \left.\textrm{ or } \{ i, l\}\neq \{d-1,d\} \textrm{ and }  j=  \left \lceil \frac{i+{l}}{2}\right \rceil +\delta_{i,{l}} \, \right\}.\end{multline*}
As for Fig. 1, 2 and 3 in \cite{aizav}, the two edges corresponding to \eqref{eq:points} are missing.

Additionally, there are some mismatches concerning Fig. 2 and $\omega_3$ in \cite{aizav}. The simplest way to make the labels of Fig. 2 match is fist to multiply $\omega_3$ in page 273 by $5$, obtaining  
\begin{equation}\label{eq:omega3}
w_3((i,j))(m)=
\begin{cases}5^{|i-j|+1} &\text{if }m\equiv i-j\textrm{ (mod 3)}\\
3 \cdot 5^{|i-j|}  &\text{if }m\equiv (i-j-sign(i-j+1/2))\textrm{ (mod 3)}\\
0 &\text{if }m\equiv (i-j+sign(i-j+1/2))\textrm{ (mod 3)},\\
\end{cases}
\end{equation}
and then consider the colouring [red is $m=0$, green is $m=2$, blue is $m=1$], so that we just have to swap the values of the blue and green labels along the diagonal of Fig. 2 in page 314. 

\begin{figure}[h!]
\centering
\includegraphics[width=\textwidth]{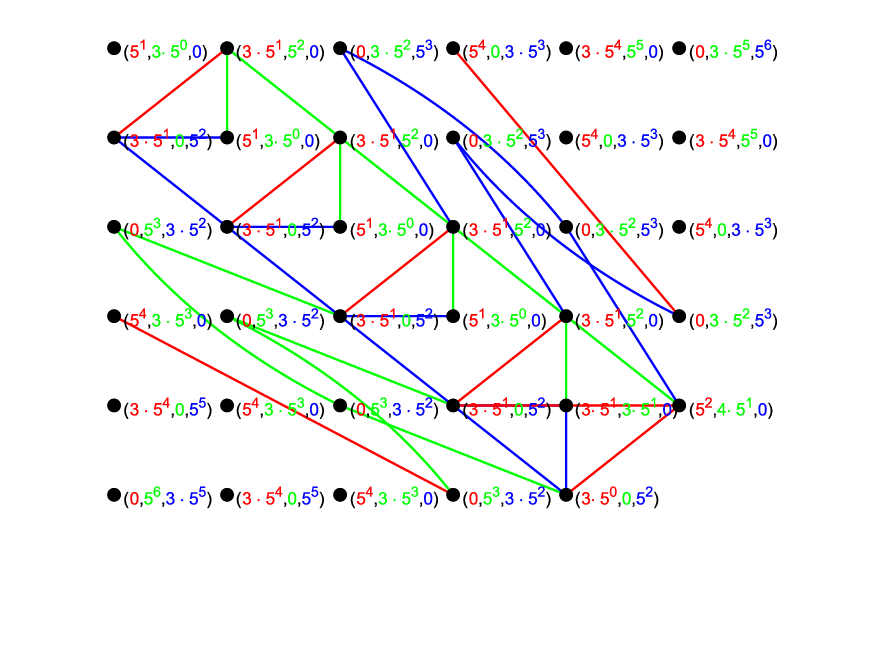}
\vspace{-2.7cm}
    \caption{Graph $\Gamma_3$ with the weights $\omega_3$ for the case $d = 6$.}
\label{fig:1}       
\end{figure}

Finally, the key to prove now Lemma 2.40 in the closest way to that of \cite{aizav} is redefining $\omega_3$ at the nodes $(d-1,d-1)$, $(d-1,d)$ and $(d,d-1)$: 
\begin{align*}
\omega_3(d-1,d-1)(0)&=3\cdot 5^1, & \omega_3(d-1,d)(0)&=5^2 ,& \omega_3(d,d-1)(0)&=3\cdot 5^0 , \\
\omega_3(d-1,d-1)(2)&=3\cdot 5^1,&\omega_3(d-1,d)(2)&=4\cdot 5^1 ,& \omega_3(d,d-1)(2)&= 0,\\
\omega_3(d-1,d-1)(1)&=0 ,& \omega_3(d-1,d)(1)&=0 ,& \omega_3(d,d-1)(1)&=5^2,
\end{align*}
so that the resulting diagram for $d=6$ (Fig. 2 in \cite{aizav}) is given by Image \ref{fig:1}.

\vspace{-0cm}

\vspace{.3cm}

With the redefinition of $\omega_3$ above, we only need to add to the proof of Lemma 2.40 in \cite{aizav} an analysis of the edges around $(d-1,d-1)$. They look, for $d\geq 3$, like the ones in Image \ref{fig:1}. We thus have forests with maximal degree $\leq 3$, as we need. Finally, the cases $d=1,2$ are straigthforward.

\newpage

\section{An alternative solution}

We indicate here how to get an alternative solution that, although requires more changes, would keep better the original intuition for the proof.

\begin{figure}[h!]
\includegraphics[width=\textwidth]{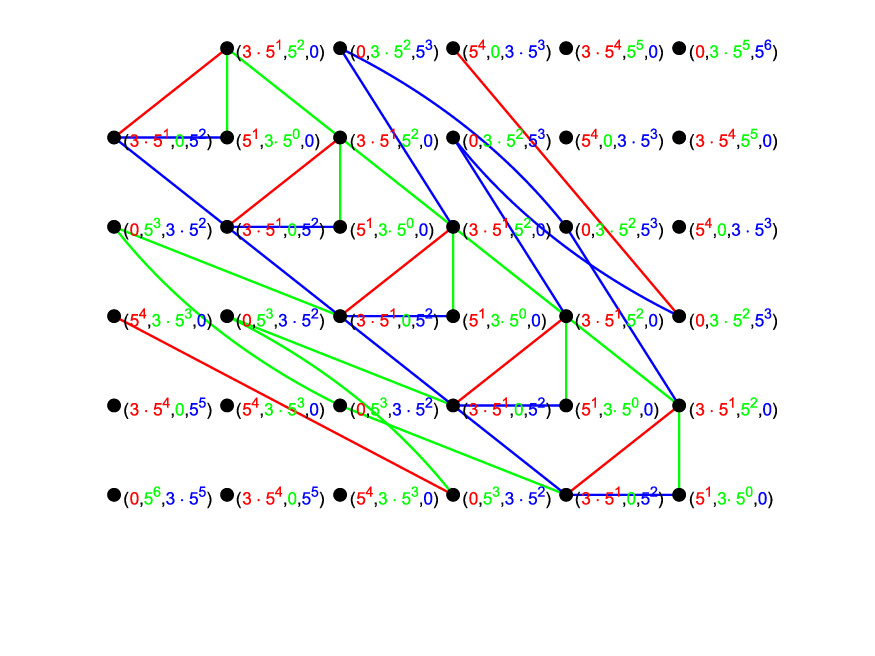}
\vspace{-2.7cm}
\caption{Graph $\Gamma_3$ with $I\neq J$ for the case $d = 6$.}
\label{fig:2}       
\end{figure}

At the beginning of Section 2.6 in \cite{aizav}, recall that $L=\{1,\ldots,d\}$. Stay with $J=L\times L\setminus \{(d,d)\}$ and replace $I$ by $I=L\times L\setminus \{(1,1)\}$. This entails  changes in $S_j$ and $\Gamma_j$, which we omit here for the sake of brevity. With the definition of $\omega_3$ as in \eqref{eq:omega3} (without any redefinition) and the same conventions for the colours as described in Section \ref{sec:correction}, the corresponding $\Gamma_3$ for $d=6$ is given by Image  \ref{fig:2}.

\vspace{-0.3cm}
\vspace{.5cm}
The general proof then follows along the same lines as the original one.

\vspace{.2cm}

\end{document}